\documentclass[preprint, 11pt]{amsart}
\usepackage{amsmath,amsthm,amsfonts,amssymb,amscd}
\usepackage{fancyhdr, tikz-cd, tikz}
\usetikzlibrary{arrows,chains,matrix,positioning,scopes}
\newtheorem{problem}{PROBLEM}
\newtheorem{theorem}{THEOREM}[section]
\theoremstyle{definition}
\newtheorem{corollary}{Corollary}[theorem]
\newtheorem{definition}[theorem]{Definition}
\newtheorem{remark}[theorem]{Remark}
\newtheorem{example}{Example}[section]
\newtheorem{lemma}[theorem]{Lemma}
\newtheorem{proposition}[theorem]{Proposition}
\numberwithin{equation}{section}
\newcommand{\C}{\mathbb{C}}

\newcommand{\K}{\mathcal{K}}

\newcommand{\E}{E=k(t_1,\cdots,t_n)}

\newcommand{\B}{\mathcal{B}}

\newcommand{\poly}{\mathcal{P}}
\newcommand{\rat}{\mathcal{R}}

\newcommand{\ok}{\overline{k}}

\renewcommand{\a}{\alpha}
\renewcommand{\b}{\beta}

\newcommand{\bt}{\begin{theorem}}
\newcommand{\et}{\end{theorem}}
\newcommand{\bco}{\begin{corollary}}
\newcommand{\eco}{\end{corollary}}
\newcommand{\bd}{\begin{definition}}
\newcommand{\ed}{\end{definition}}
\newcommand{\bp}{\begin{problem}}
\newcommand{\ep}{\end{problem}}
\newcommand{\bl}{\begin{lemma}}
\newcommand{\el}{\end{lemma}}
\newcommand{\bprop}{\begin{proposition}}
\newcommand{\eprop}{\end{proposition}}
\newcommand{\br}{\begin{remark}}
\newcommand{\er}{\end{remark}}
\newcommand{\bpf}{\begin{proof}}
\newcommand{\epf}{\end{proof}}
\newcommand{\bex}{\begin{example}}
\newcommand{\eex}{\end{example}}
\title[liouvillian solutions of differential equations]{liouvillian solutions of first order non linear differential equations}
\author{Varadharaj Ravi  Srinivasan}

\address{Department of Mathematical Sciences \\
IISER Mohali, SAS Nagar, Punjab 140306}

\begin{document}

\maketitle

\begin{abstract}

Let $k$ be a differential field of characteristic zero and $E$ be a liouvillian extension of $k$. For any differential subfield $K$ intermediate to $E$ and $k$, we prove that there is an element in the set $K-k$ satisfying a linear homogeneous differential equation over $k$. We apply our results to study liouvillian  solutions of first order non linear differential equations and provide generalisations and new proofs for several results of M. Singer and M. Rosenlicht on this topic.

\end{abstract}


\section{Introduction}  

Throught this article, we fix  a differential field $k$ of characteristic zero.  Let $E$ be a differential field extension of $k$ and let $'$ denote the derivation on $E$. We say that  $E$ is a \textit{ liouvillian extension} of $k$  if $\E$ and there is a tower of differential fields 
\begin{equation*}k=k_0\subset k_1\subset\cdots\subset k_n=E\end{equation*}
such that for each $i$,  $k_i=k_{i-1}(t_i)$ and either $t'_i\in k_{i-1}$ or $t'_i/t_i\in k_{i-1}$ or $t_i$ is 
algebraic over $k_{i-1}$. If $\E$ is a liouvillian extension of $k$ such that $t'_i\in k_{i-1}$ for each $i$ then we call $E$ an \textit{iterated antiderivative extension} of $k$. A solution of a differential equation over $k$  is said to be  liouvillian over $k$ if the solution belongs to  some liouvillian extension of $k$.  

Let $\poly$ be an $n+1$ variable polynomial over $k$. We are concerned with the  liouvillian solutions of the differential equation \begin{equation}\label{firstdiffereqn}\poly(y,y',\cdots, y^{(n)})=0.\end{equation}   We prove in theorem \ref{solutiondense} that \begin{itemize}
\item[] if $E$ is a liouvillian extension of $k$ and  $K$ is a differential field intermediate to $E$ and $k$ then $K=k\langle u_1,\cdots,u_l\rangle$ where for each $i$, the element $u_i$ satisfies a linear  differential equation over $k\langle u_1,\cdots, u_{i-1}\rangle$. Moreover if $E$ is an iterated antiderivative extension of $k$ having the same field of constants as $k$ then each $u_i$ can be chosen so that $u'_i\in k(u_1,\cdots, u_{i-1})$, that is,   $K$ is also an iterated antiderivative extension of $k$. \end{itemize}

Our result regarding iterated antiderivative extensions generalises the main result of \cite{Sri2010} to differential fields $k$ having a non algebraically closed  field of constants. The main ingredient used in the proof our theorem is the lemma \ref{ldesol} using which we also obtain the following interesting results concerning solutions of non linear differential equations.  
\begin{itemize}

\item[A.] In remark \ref{order2}, we show that if $E$ is an iterated antiderivative extension of $k$ having the same field of constants and if $y\in E$ and $y\notin k$ satisfies a differential equation $y'=\poly(y)$, where $\poly$ is a polynomial in one variable over $k$, then degree of $\poly$ must be less than or equal to $2$. \\
\item[B.] Let $C$ be an algebraically closed field of characteristic zero with the trivial derivation.  In proposition \ref{M-B}, we prove that for any rational function $\rat$ in one variable over $C$, the differential equation $y'=\rat(y)$ has a non constant liouvillian solution $y$  if and only if $1/\rat(y)$ is of the form $\partial z/\partial y$ or $(1/az)(\partial z/ \partial y)$ for some $z\in C(y)$ and for some  non zero element $a\in C$.   This result generalises and provides a new proof for a result of Singer (see \cite{Sin1975}, corollary 2). We also prove that 
 for any polynomial $\poly$ in one variable over $C$ such that the degree of $\poly$ is greater than or equal to $3$ and that $\poly$ has no repeated roots,   the differential equation $(y')^2=\poly(y)$ has no non constant liouvillian solution over $C$. This result appears as proposition \ref{hyperellipticcurves} and it generalises an observation made by Rosenlicht \cite{M.Ros} concerning non constant liouvillian solutions of the elliptic equation $(y')^2=y^3+ay+b$ over complex numbers with a non zero discriminant. \\
 \item[C.]  Using theorem \ref{solutiondense} one can construct a  family of differential equations with only algebraic solutions: Let $\a_2,\a_3,$ $\cdots,\a_n\in k$ such that $x'\neq \a_2$ and $x'\neq \a_3$ for any $x\in k$. Let $\ok$ be  an algebraic closure of $k$ and let  $E$ be a  liouvillian extension  of $\ok$ with $C_E=C_{\ok}$. We prove in proposition \ref{algsolutions} that if there is an element $y\in E$  such that  $$y'=\a_ny^n+\cdots+\a_3y^3+\a_2y^2$$ then    $y\in \overline{k}$.   

\end{itemize}

In a future publication, the author hopes to develop the techniques in this paper further to provide an algorithm to solve the following problem, which appears as ``Problem 7" in \cite{Sin1990}: Give a procedure to decide if a polynomial first order differential equation $\poly(y,y')=0$, over the ordinary differential field $\C(x)$ with the usual derivation $d/dx$,  has an elementary solution and to find one if it does.

{\bf Preliminaries and Notations.}
 A \textsl{derivation} of the field $k$, denoted by $'$, is an additive endomorphism of $k$ that satisfies the Leibniz law $(xy)'=x'y+xy'$ for every $x,y\in k$. A field equipped
with a derivation map is called a \textsl{differential field}. For any $y\in k$, we will denote first and second derivative of $y$ by $y'$ and $y''$ respectively and for $n\geq 3$, the $n$th derivative of $y$ will be denoted by $y^{(n)}$. The set of \textsl{constants} $C_E$ of a differential field $E$ is the kernel of the endomorphism $'$ and it can be seen that the set of constants is a differential subfield of $k$.  Let $E$ and $k$  be differential fields. We say that $E$ is a \textsl{differential field extension}  of $k$ if $E$ is a field extension of $k$ and the restriction of the derivation of $E$ to $k$ coincides with the derivation of $k$. Whenever we write $k\subset E$ be differential fields, we mean that $E$ is a differential field extension of $k$ and we write $y\in E-k$ to mean that $y\in E$ and  $y\notin k$. The transcendence degree of a field extension $E$ of $k$ will be denoted by tr.d$(E|k)$. If $E$ is a field (respectively a differential field),  $M$ is a subfield (respectively a differential subfield) of $E$ and $K$ is a subset of $E$  then  the smallest subfield (respectively the smallest differential subfield) of $E$ containing both  $M$ and $K$ will be denoted by $M(K)$ (respectively by $M\langle K\rangle$). It is easy to see that the field $M(K)$ is a differential field if both $M$ and $K$ are differential subfields of a differential field.   It is well known that every derivation of $k$ can be uniquely extended to a derivation of any algebraic extension of $k$ and in particular, to any algebraic closure $\overline{k}$ of $k$. Thus if $E$ is a differential field extension of $k$ and tr.d$(E|k)=1$ then for any $y\in E$ transcendental over $k$ the derivation $\partial/\partial y$ on $k(y)$  uniquely extends to a derivation of  $E$.   We refer the reader to \cite{Kap}, \cite{Mag} and \cite{Put-Sin} for basic theory of differential fields and for the reader's convenience and for easy reference, we record a basic result concerning liouvillian extensions  in the following theorem.
\bt\label{construction} Let $k$ be a differential field,  $E$ be a field extension of $k$ and $w\in E$ be transcendental over $k$. For an element $\a\in k$ if there is no element $x\in k$ such that $x'=\a$ then there is a unique derivation on $k(w)$ such that $w'=\a$ and that $C_{k(w)}=C_k$.  Similarly, if there is no element $x\in k$ such that $x'=l\a x$ for any positive integer $l$ then there is a unique derivation on $k(w)$ such that $w'=\a w$ and that $C_{k(w)}=C_k$. Finally, if $E$ is a differential field extension of $k$ having an element $u$ algebraic over $k$ such that $u'\in k$ then there is an element $\b\in k$ such that $u'=\b'$.   \et


\section{Role of First Order Equations}

In this section we prove our main result.  The proof is based on the theory of linearly disjoint fields  and in particular, we will heavily rely on the following  lemma. 


\lemma \label{ldesol}
Let $k\subset E$ be differential fields and let $K$ and $M$ be  differential fields  intermediate to  $E$ and $k$ such that $K$ and $M$ are linearly disjoint over $k$ as fields. Suppose that there is an element $y\in M(K)-M$ such that $y'=fy+g$ for some $f,g\in M$. Then there exist a monic linear  differential polynomial $L(Y)$ over $k$  of degree $\geq 1$ and an element $u\in K-k$ such that $L(u)=0$.

\bpf

Let $\B:=\{e_\a\ |\ \a\in J\}$ be a  basis of the $k-$vector space $M$. Since $M$ and $K$ are linearly disjoint,  $\B$ is also a basis of $M(K)=K(M)$ as a $K-$vector space.   There exist elements  $\a_{1},\cdots,\a_r,$ $\a_{r+1},\cdots,\a_t\in J$ such that  
\begin{equation}
y=\sum^r_{j=1}u_je_{\a_j},
\end{equation}
where $u_j\in K$ and for $i=1,\cdots, r,$ 
\begin{equation}
e'_{\a_i}=\sum^t_{p=1}n_{pi}e_{\a_p},\quad fe_{\a_i}=\sum^t_{p=1}l_{pi}e_{\a_p},\quad g=\sum^t_{p=1}m_pe_{\a_p},
\end{equation}
 where each $n_{pi}, l_{pi}$ and $m_p$ belongs to $k$.

Then \begin{align*}
y'&=\sum^r_{j=1}u'_je_{\a_j}+\sum^r_{j=1}u_j\sum^t_{p=1}n_{pj}e_{\a_p}\\
&=\sum^r_{j=1}u'_je_{\a_j}+\sum^t_{p=1}\left(\sum^r_{j=1}n_{pj}u_j\right)e_{\a_p}\\
&=\sum^r_{p=1}\left(u'_p+\sum^r_{j=1}n_{pj}u_j\right)e_{\a_p}+\sum^t_{p=r+1}\left(\sum^r_{j=1}n_{pj}u_j\right)e_{\a_p}
\end{align*}
and $fy=\sum^r_{j=1}u_jfe_{\a_j}$ $=\sum^t_{p=1}\left(\sum^r_{j=1}l_{pj}u_j\right)e_{\a_p}$.

Now, from the equation $y'=fy+g$, we obtain for $p=1,\cdots,r,$ that
\begin{equation}
u'_p+\sum^r_{j=1}n_{pj}u_j=\sum^r_{j=1}l_{pj}u_j+m_p.
\end{equation}
Consider the $k-$vector space $\K:=span_k\{1,u_1,\cdots,u_r\}$. From the above equation, it is clear that $\K$ is a (finite dimensional) differential $k-$vector space and that $k\subset \K\subset K$.  Then for every $u\in \K-k$, there exists a non negative integer $n$ such that  $u=u^{(0)}, u^{(1)},\cdots,u^{(n)}$  are linearly independent over  $k$ and that $u^{(n+1)}=\sum^n_{i=0}a_iu^{(i)}$, for some  $a_i\in k$.  Hence the lemma is proved.
\epf



\bt\label{solutiondense}
Let $E$ be a liouvillian extension field of $k$  and let $K$ be a differential field intermediate to $E$ and $k$.  Then \begin{itemize} \item [I.]$K=k\langle u_1,\cdots,u_l\rangle$ where for each $i$, $u_i$ satisfies a linear homogeneous differential equation over $k\langle u_1,\cdots, u_{i-1}\rangle$.\\ \item[II.] If $\E$ is an iterated antiderivative extension of $k$ with $C_E=C_k$  then $K$ is an iterated antiderivative extension of $k$. That is, $K=k(v_1,\cdots,v_m)$, where $v'_i\in k(v_1,\cdots,v_{i-1})$ .\\   \item[III.] If $k$ is an algebraically closed field  then there is  an element $u\in K-k$ such that $u'=au+b$ for some $a,b\in k$.\\ \item[IV.]Finally, if $k$ is algebraically closed field and $c'=0$ for all $c\in k$,  then there is an element $z\in K-k$ such that either $z'=1$ or $z'=az$ for some element $a\in k$. \end{itemize}
\et

\bpf 
We will use an induction on tr.d$(E|k)$ to prove item I.
For any differential field $k^*$ intermediate to $E$ and $k$ we observe that $E$ is a liouvillian extension of $k^*$ and therefore for differential fields $E\supset K\supset k^*$  such that tr.d$(E|k^*)<$tr.d$(E|k)$, we shall assume that the item I of our theorem holds.  
If $u\in K$ is algebraic over $k$ then the vector space dimension $[k(u):k]=m$ for some positive integer $m$. Since $k(u)$ is a differential field,  the set $\{u, u',u'',\cdots, u^{(m)}\}$ is linearly dependent and therefore $u$ satisfies a linear differential equation over $k$. Since
 $k$ is of characteristic zero, there exist an element  $u_1\in K$ such that that $k(u_1)$ is the algebraic closure of $k$ in $K$.  Replacing $k$ with $k(u_1)$, if necessary, we shall assume that $k$ is algebraically closed in $K$ and that tr.d$(K|k)\geq 1$. 
 Let $\E$, $k_0:=k$ and choose the largest positive integer $m$ so that $k_{m-1}$ and $K$ are linearly disjoint over $k$. Since $k$ is  algebraically closed in $K$ and the characteristic of $k$ is zero, the field $k_{m-1}$ is algebraically closed in  $k_{m-1}(K)$. Now our choice of $m$ guarantees that $t_m$ is not algebraic over $k_{m-1}$ and that $t_m$ is algebraic over $k_{m-1}(K)$.   Let $\poly(X)=\sum^n_{i=0}a_{n}X^n$ be the monic irreducible polynomial of $t_m$ over $k_{m-1}(K)$ and let $l$ be the smallest integer such that $a_{n-l}\notin k_{m-1}$. Expanding the equation $(\sum^n_{i=0}a_{i}t^i_m)'=0$  and comparing it with the equation $\sum^n_{i=0}a_{i}t^i_m=0$, we obtain   
 
 $$a'_{n-l}= \begin{cases} -(n-(l-1))ra_{n-(l-1)} &\ \mbox{if }\  t'_m=r\in k_{m-1} \\ 
(n-l)ra_{n-l} &\ \mbox{if }\ t'_m/t_m=r\in k_{m-1}. \end{cases} $$
  
Thus there is an element $a_{n-l}\in k_{m-1}(K)-k_{m-1}$ such that  either $a'_{n-l} \in k_{m-1}$ or $a'_{n-l}/a_{n-l}\in k_{m-1}$. Now we apply lemma \ref{ldesol} and obtain an element $u_2\in K-k$ such that $L(u_2)=0$, where $L(Y)=0 $ is a linear  differential equation over $k$ of order $\geq 1$. Since $k$ is algebraically closed in $K$, such an element $u_2$ must be transcendental over $k$ and we now prove item I by setting $k^*=k\langle u_2\rangle$ and invoking the induction hypothesis.
 
Let $E$ be an iterated antiderivative extension of $k$ with $C_E=C_k$. We shall choose elements $t_1,\cdots,t_n\in E$ so that $\E$ and that $t_1,\cdots,t_n$ are algebraically independent over $k$ (see \cite{Sri2010}, theorem 2.1).  Let $L(Y)=0$ be a linear homogeneous differential equation over $k$ of smallest positive degree $n$ such that $L(u)=0$ for some $u\in K-k$. Let $u\in k_{i}-k_{i-1}$ and observe that the ring $k_{i-1}[t_i]$ is a differential ring with no non trivial differential ideals (see \cite{Put-Sin}, example 1.18). Consider the differential $k-$vectorspace $W:=$span$_k\{u,u',\cdots,u^{(n-1)}\}$ and choose a nonzero  element $g\in k_{i-1}[t_i]$ so that $gu^{(j)}\in k_{i-1}[t_i]$  for all $j$. Then  the set $I=\{h\in k_{i-1}[t_i]\ |\ hW\subset k_{i-1}[t_i]\}$ is a non zero differential ideal of $k_{i-1}[t_i]$ and therefore $I=R$. Thus $1.u=u\in k_{i-1}[t_i]$ and let $u=a_pt^p_i+\cdots+a_0$, where $a_p\neq 0$ and for each $j$, $a_j\in k_{i-1}$. Then $L(a_pt^p_{i}+\cdots+a_0)=0$ together with the fact that $t_i$ is transcendental over $k_{i-1}$ implies $L(a_p)=0$ and thus we have found a non zero solution $a_p$ of $L(Y)=0$ in $k_{i-1}$. Repeating this argument, one can then show that there is a non zero element $d\in k$ such that $L(d)=0$. This means that $L(Y)=L_{n-1}(L_1(Y))$, where $L_1(Y)=Y'-(d'/d)Y$ and $L_{n-1}(Y)=0$ is a linear  homogeneous differential equation over $k$ of order $n-1$. Since $L_1(u)\in K$, from our choice of $u$ and $n$, we obtain that $L_1(u)\in k$. Let $v=u/d$ and observe that $v\in K-k$ and that $v'\in k$. Since $u\in k_i-k_{i-1}$, we know that $t_i$ is algebraic over $k_{i-1}(v)$ and since $C_E=C_k$, it follows from theorem \ref{construction} that  $t_i\in k_{i-1}(v)$. Thus $k(v)\subset K\subset E=k(v)(t_1,\cdots,t_{i-1},t_{i+1},\cdots,t_n)$ and an induction on the transcendence degree, as in the proof of I, will prove II.

Suppose that the differential field $k$ is an algebraically closed field. Let $L(Y)=0$ be a linear homogeneous differential equation over $k$ of smallest positive degree $n$ such that $L(u)=0$ for some $u\in K-k$. Then $L(Y)=0$ admits a liouvillian solution and  since $k$ is algebraically closed, the field of constants $C_k$ is algebraically closed as well. It follows from \cite{Sin1981}, theorem $2.4$ that  $L(Y)=L_{n-1}(L_1(Y))$, where  $L_{n-1}(Y)=0$ and $L_1(Y)=0$ are linear homogeneous differential equations over $k$ of degrees $n-1$ and $1$ respectively. Thus from the choice of $L$, we must have $L_1(u)\in k$.

Suppose that $k$ is an algebraically closed field with the trivial derivation. Assume that there is no $z\in K-k$ such that $z'=0$. We know from II that there is an element $u\in K-k$ such that $L(u)=0$, where $L(Y)=Y'-bY-c$ for some $b,c\in k$. If $b=0$ then let $z=(u/c)$ and observe that $z'=1$. Therefore we shall assume that $b\neq 0$. We will now show that there is an element $z\in k(u)-k$ such that $z'=mb z$ for some integer $m>0$.  Suppose  that there is no such integer $m$. Then we shall consider the field $E=k(u)(X)$, where $X$ is transcendental over $k(u)$ and define a derivation on $E$ with $X'=bX$. Then from  theorem \ref{construction} we obtain that  $C_E=C_{k(u)}$. But 
$$\left(\frac{u}{X}\right)'=\left(\frac{-c}{bX}\right)'$$
and therefore $(bu+c)(1/bX)\in C_{k(u)}\subset k(u)$. Thus we obtain $X\in k(u)$, which is absurd. Hence the theorem is proved. \epf


\br\label{order2}

Let $E$ be an iterated antiderivative extension of $k$ with $C_E=C_k$. Then $\E$, where $t_1,\cdots,t_n$ are algebraically independent over $k$ and furthermore, the differential field $E$ remains a purely transcendental extension over any of its differential subfields containing $k$ (see \cite{Sri2010}, theorems 2.1 \& 2.2). Suppose that  $y\in E-k$ and that $y'=\poly(y)$ for some polynomial $\poly$ in one variable over the differential field $k$. We will now show that deg $\poly\leq 2$. The field $k(y)$ is a differential field intermediate to $E$ and $k$ and therefore, by theorem \ref{solutiondense}, there is an element $z\in k(y)-k$ such that $z'\in k$. Moreover, the element $y$ is algebraic over the differential subfield $k(z)$. Then since $E$ is a purely transcendental extension of $k(z)$, we obtain  that $y\in k(z)$. Thus $k(y)=k(z)$ and from  L{\"u}roth's theorem,  we know that there are elements $a,b,c,d\in k$ such that $ad-bc\neq 0$ and that $$z=\frac{ay+b}{cy+d}.$$ 
Let $\poly(y)= a_ny^n+a_{n-1}y^{n-1}+\cdots+a_0$, $a_n\neq 0$ and $z'=r\in k$ and observe that \begin{equation*}(a'y+b'+a\poly(y))(cy+d)-(c'y+d'+c\poly(y))(ay+b)=r(cy+d)^2.\end{equation*}
Since $a_n\neq 0$ and $ad-bc\neq 0$, comparing the coefficients of $y^n$ in the above equation, we obtain that $n\leq 2$.  On the other hand, the differential equation $Y'=-rY^2$ has a  solution in $E-k$, namely, $1/z$ and thus  the bound $n\leq 2$ is sharp.

\er

\section{First Order  Non Linear Differential Equations} 

Let $C$ be  an algebraically closed  field of characteristic zero and view $C$ as a differential field with the trivial derivation $c'=0$ for all $c\in C$.  Consider the field $C(X)$ of rational functions in one variable $X$. We are interested in the non constant liouvillian solutions of the following differential equations  \begin{enumerate} \item[(i)] $y'=\rat(y)$, where $\rat(X)\in C(X)$ is a non zero element and \item[(ii)] $(y')^2=\poly(y),$ where $\poly(X)\in C[X]$ and $\poly(X)$ has no repeated roots in $C$.  \end{enumerate}

Suppose that there is a non constant element $y$ in some differential field extension  of  $C$ satisfying either (i) or (ii) then $C(y,y')$ must be a differential field, the element $y$ must be transcendental over $C$ and for any $z\in C(y)$, we have 
\begin{equation}\label{two derivations} z'=y' \frac{\partial z}{\partial y}. \end{equation}

In the next proposition we will generalise a result of Singer (see \cite{Sin1975}, corollary 2) concerning elementary solutions of first order differential equations. 
\bprop\label{M-B}
Let $C$ be  an algebraically closed  field of characteristic zero with the trivial derivation and let $\rat(X)\in C(X)$ be a non zero element.  The equation $Y'=\rat(Y)$  has a non constant solution $y$ which is liouvillian over $C$ if and only if there is an element $z\in C(y)$ such that
$$\dfrac{1}{\rat(y)}\quad\text{is of the form}\quad \frac{\partial z}{\partial y} \ \ \text{or}\ \
										\dfrac{\frac{\partial z}{\partial y}}{a z}$$
for some non zero element $a\in C$.

\eprop

\bpf
 Let $y$ be a non constant liouvillian solution of $Y'=\rat (Y)$. From equation \ref{two derivations}, we observe that if $z\in C(y)$ and $z'=0$ then $\partial z/\partial y=0$ and since the field of constants of $C(y)$ with the derivation $\partial/\partial y$ equals $C$, we obtain that $z\in C$. Since $y$ is liouvillian, the differential field $M=C(y)$ is contained in some liouvillian extension field of $C$ and therefore, applying  theorem \ref{solutiondense}, we  obtain an element $z\in C(y)-C$ such that either $z'=1$ or $z'=az$ for some (non zero) $a\in C$. Now it follows immediately from equation \ref{two derivations} that $1/\rat(y)$ has the desired form.  Let $y'=\rat(y)$ and $1/\rat(y)$ equals $\partial z/\partial y$ or $(1/az)(\partial z/ \partial y)$ for some $z\in C(y)$ and for some non zero element $a\in C$. Then since $z'= y'(\partial z/\partial y)$,  we obtain $z'=1$ or $z'=az$.  From the fact that $C(y)\supset C(z)$, we see that $y$ is algebraic over the liouvillian extension $C(z)$ of $C$ and thus $C(y)$ is a liouvillan extension of $C$.  \epf


\bprop\label{hyperellipticcurves}

Let $C$ be  an algebraically closed  field of characteristic zero with the trivial derivation. Let $\poly(X)\in C[X]$ be a polynomial of degree $\geq 3$  with no repeated roots.  Then the differential equation $(Y')^2=\poly(Y)$ has no non constant liouvillian solution over $C$. In particular, the elliptic function $y$ such that $(y')^2=y^3+ay+b$, where $\frac{a^3}{27}+\frac{b^2}{4}\neq 0$, is not  liouvillian.

\eprop

\bpf
Suppose that there is a non constant liouvillian solution $y$ satisfying   the equation $(Y')^2=\poly(Y)$. Since $C$ is algebraically closed, such an element $y$ must be transcendental over $C$.  Applying theorem \ref{solutiondense} to the differential field $C(y,y')$, we obtain that there is an element $z\in C(y,y')-C$ such that either $z'=1$ or $z'=az$ for some  $a\in C$. We will first show that there is no $z\in C(y,y')-C$ with $z'=az$ for any $a\in C$. 

For ease of notation, let $P=\poly(y)$  and observe that $P$ is not a square in $C[y]$ and thus $y'\notin C(y)$ and therefore $y'$ lies in a quadratic extension of $C(y)$.
Write $z=A+By'$, where $A,B\in C[y]$ and $B\neq 0$. Then, taking derivatives, we obtain $A'+B'y'+By''=z'$ and using equation \ref{two derivations}, we  obtain  $$y' \frac{\partial A}{\partial y} + P\frac{\partial B}{\partial y}+  \frac{B}{2} \frac{\partial P}{\partial y}=z'.$$
If there is a non zero element $z\in C(y,y')$ such that $z'=az$ then, by comparing coefficients of the above equation, we obtain \begin{align}&\frac{\partial A}{\partial y}=aB\label{first eqn}\qquad\text{and}\\  &P\frac{\partial B}{\partial y}+  \frac{B}{2} \frac{\partial P}{\partial y}=aA\label{second eqn}.\end{align} 
Multiplying the equation \ref{second eqn} by $2B$ and using \ref{first eqn}, we obtain $\partial (B^2 P)/\partial y=\partial A^2/\partial y. $
Thus there is a non zero constant $c\in C$ such that \begin{equation}\label{hypeleqn}B^2P=A^2+c.\end{equation} Write $A=A_1/A_2$ and $B=B_1/B_2$, where $A_1,A_2,B_1,B_2$ are polynomials in $C[y]$ such that  $A_1,A_2$ are relatively prime, $B_1,B_2$ are relatively prime  and $A_2,B_2$ are monic. Then since $P$ has no square factors, it follows from the equation $$B^2_1A^2_2P=A^2_1B^2_2+cA^2_2B^2_2$$ that $A_2=B_2$. The equation \ref{first eqn}, together with our assumption that $A_2$ and $B_2$ are monic, forces $A_2=B_2=1$ and we obtain deg $A $= 1+deg $B$. Now from equation \ref{hypeleqn}, we have $$2\ \text{deg}\ A -2+\ \text{deg}\  P=2\ \text{deg} \ A$$ and thus we obtain deg $ P=2$, which contradicts our assumption. 

Now let us suppose that $z\in C(y,y')-C$ such that $z'=1$. Then  $\partial A/\partial y=0$, which imples $A\in C$ and we have $(\partial B/\partial y)P+  (B/2) (\partial P/\partial y)=1$. Thus  \begin{equation}\label{elliptic} \frac{\partial (B^2P)}{\partial y}=2B.\end{equation}
It is clear from the above equation that $B$ cannot be a  polynomial and that $B$ cannot have a pole of order $1$.   On the other hand if $c\in C$ is a pole of $B$ of order $m\geq 1$ then we shall write $B=R+\sum^m_{i=1}\b_i/(x-c)^i$, where $R\in C(y)$ and $c$ is not a pole of $R$. Since $P$ has no repeated roots, we conclude that $B^2P$ has  a pole at $c$ of order $\geq 2m-1$ and thus $c$ is a pole of $\partial (B^2P)/\partial y$  of order $\geq 2m$. This contradicts equation \ref{elliptic}.
\epf


 Consider the ring of rational functions in two variables $C[Y,X]$. For constants $a\in C-\{0\}$ and $b\in C$, define a derivation $D$ of $C[Y,X]$ by setting $D(Y)=X$ and $D(X)=a/2$. Then $D(X^2-aY-b)=0$ and therefore the ideal $I=\langle X^2-(aY+b)\rangle$ is a differential ideal as well as a prime ideal. Thus the factor ring $C[Y,X]/I$ is a differential ring as well as a domain. Extend the derivation to the field of fractions $E$ of $C[Y,X]/I$.  Now, in the differential field $E$, we have  elements $x,y$ such that $y'=x$, $x'=a/2$ and $(y')^2=ay+b$. Let $z=2x/a$ and note that $z'=1$. Since tr.d$(E|C)=1$, the field $E$ is algebraic over $C(z)$ and thus $E$ is liouvillian over $C$. Thus the equation $(Y')^2=aY+b$ admits non constant liouvillian solutions. For the polynomial $f(X)=1-X^2\in \C(X)$, where $\C$ is the field of complex numbers, we have $(\frac{d}{dx}(\sin x))^2=f(\sin x)$ and thus $n\geq 3$ is necessary for the proposition \ref{hyperellipticcurves} to hold.



In the next proposition, we shall construct differential equations  whose liouvillian solutions, from liouvillian extensions having the same field of constants, are all algebraic over the ground field. 

\bprop\label{algsolutions}
Let $\ok$ be  an algebraic closure of $k$ and let $E$ be a  liouvillian extension of $\ok$ with $C_E=C_{\ok}$. Let $$F(Y,Y')=Y'-\a_nY^n-\cdots-\a_{2}Y^{n-1}-\a_1Y,$$
where $\a_1,\a_2,$ $\cdots,\a_n\in k$ and suppose that there is  an element $y\in E-\ok$ such that  $F(y,y')=0$.   
\begin{itemize}
\item [I.] If there is an element $\gamma\in \ok$ such that $\gamma'=\a_1 \gamma$ then there is no $z\in \ok(y)-\ok$ such that $z'/z\in \ok$.
\item[II.] If $\a_1=0$ and $x'\neq \a_2$ for all $x\in k$ then there is an element $w\in \ok(y)-\ok$ and an element $v\in \ok$ such that $w'=\a_2$ and that $v'=\a_3$.

\end{itemize}

\eprop

\bpf
Suppose that there is an element $z\in\ok(y)-\ok$ such that $z'=\a z$ for some $\a\in \ok$. Write $z=P/Q$ for some relatively prime polynomials $P=\sum^{s}_{i=0}a_iy^i$ and $Q=\sum^t_{j=0}b_jy^j$ in $\ok[y]$ with $b_t=1$. Taking derivative, we obtain $\a PQ=P'Q-PQ'.$ Replacing $z$ and $\a$ by $1/z$ and $-\a$,  if necessary,  we shall assume that $b_0\neq 0$. Let $r$ be the smallest non zero integer such that $a_r\neq 0$.  Then we have  \begin{align}&\a(a_ry^{r}+ \cdots ) (b_0+\cdots)= (b'_0+\cdots)(a_ry^r+\cdots)\\&-((a'_r+ra_r\a_1)y^r+\cdots)(b_0+\cdots) \notag
\end{align}
 Now we compare the coefficients of $y^r$ and obtain that  $$b'_0a_r-a'_rb_0-\a a_rb_0=r\a_1 a_r b_0$$ and thus  $\left(\gamma^r a_r z/b_0\right)'=0.$ This contradicts our assumption that $C_{E}\neq C_{\ok}$ and thus item I is proved. 
 
Let us assume that $\a_1=0$ and that $x'\neq \a_2$ for all $x\in k$. From theorem \ref{solutiondense}, we obtain an element $z\in \ok(y)-\ok$ with $z'=\a z+\b$ for some $\a,\b\in \ok$. We claim that there is an element $z_1\in \ok(y)-\ok$ with $z'_1\in \ok.$ Write $z=P/Q$ for some relatively prime polynomials $P=\sum^{s}_{i=0}a_iy^i$ and $Q=\sum^t_{j=0}b_jy^j$ in $\ok[y]$ with $b_t=1$. Taking derivative, we obtain \begin{equation}\label{Abel} \a PQ+\b Q^2=P'Q-PQ'.\end{equation}
Comparing the constant terms of the above equation, we have $\a a_0b_0+\b b^2_0=a'_0b_0-a_0b'_0$. If $b_0\neq 0$ then  $(a_0/b_0)'=\a(a_0/b_0)+\b$ and thus $[z-(a_0/b_0)]'=\a[z-(a_0/b_0)]$. This contradicts item I of the proposition.   Thus $b_0=0$ and  we choose the smallest integer $m$ such that $b_m\neq 0$ and observe that  $a_0\neq 0$. From equation \ref{Abel}, we have 
\begin{align}\label{inhomo}&\a(b_my^{m}+ b_{m+1}y^{m+1}+\cdots ) (a_0+a_1y+\cdots)+\b(b^2_my^{2m}+ \\ & 2b_mb_{m+1}y^{2m+1}+\cdots ) = (a'_0+a'_1y+\cdots)(b_my^m+\cdots)
-(b'_my^m +\notag\\ &(mb_m\a_2+b'_{m+1})y^{m+1}+\cdots)(a_0+a_1y+\cdots) \notag
\end{align}
Compare the coefficients of $y^m$ and obtain $\a a_0b_m=a'_0b_m-a_0b'_m$.   Therefore $(a_0/b_m)'=\a (a_0/b_m)$ and since   $z'=\a z+\b$, it follows that
\begin{equation} \label{antipresent}\left(\frac{b_mz}{a_0}\right)'= \frac{\b b_m}{a_0}.\end{equation}
Taking $z_1=b_mz/a_0$, we prove our claim.  Thus, in the above calculations,  we shall suppose that $\a=0$. Then we have \begin{align}\label{antieqn}& \b(b^2_my^{2m}+ 2b_mb_{m+1}y^{2m+1}+\cdots )=(a'_0+ a'_1y+\cdots)(b_my^m+  b_{m+1}y^{m+1}\\ &+\cdots)-  (a_0+a_1 y+\cdots) (b'_my^m+(mb_m\a_2+b'_{m+1})y^{m+1}+\cdots) \notag
\end{align}
Equating the coefficients of $y^m$, we obtain  $a'_0b_m-b'_ma_0=0$ and thus  for some non zero  $c\in C_{\ok}$, we have $ca_0=b_m$. Note that if $m\geq 2$ then comparing the  coefficient  of $y^{m+1}$, we obtain \begin{equation}\label{inhomo2}a'_0b_{m+1}-a_0b'_{m+1}+a'_1b_m-a_1b'_m=mb_m\a_2a_0\end{equation} and substituting $b_m=ca_0$, we obtain that $[(b_{m+1}+ca_1)/(mca_0)]'=\a_2$. Since $(b_{m+1}+ca_1)/(mca_0)\in \ok$, we shall apply theorem \ref{construction} and obtain a  contradiction to our assumption on $\a_2$. Therefore $m=1$ and from equation \ref{antieqn}, we have $f'=\a_2+c\b$, where $f=(b_{2}+ca_1)/(ca_0)$. Then $(f-cz)'=\a_2$, where $f-cz\in \ok(y)$ and $f-cz\notin \ok$.  This proves the first part of  item II.

Now, in all the above calculations, we shall assume that $\a=0$ and $\b=\a_2$. Then $f'=(1+c)\a_2$ and note that if $c\neq -1$ then $(f/(1+c))'=\a_2$ and again we obtain a contradiction to our assumption on $\a_2$. Thus $c=-1$.  Now, $f'=0$ and consequently we have $a_1-b_2=c_1a_0$ for some constant $c_1\in C_{\ok}$. We compare the coefficients of $y^3$ and obtain \begin{align}a'_2b_1-a_2b'_1+a'_1b_2-a_1b'_2+a'_0b_3-a_0b'_3= 2b_1b_2\a_2+2a_0b_2\a_2+a_0b_1\a_3\notag.\end{align} Substituting $b_1=-a_0$ and $b_2=a_1-c_1a_0$ in the above equation, we obtain
\begin{align}-a'_2a_0+a_2a'_0-c_1a'_1a_0+c_1a_1a'_0+a'_0b_3-a_0b'_3= -a^2_0\a_3\notag.\end{align}
Then it is easy to see that  $v'=\a_3$ for $v=(a_2+c_1a_1+b_3)/a_0\in \ok$.  \epf

\subsection*{Final Comments}  One can slightly extend the proposition \ref{algsolutions} to include differential equations of the form  $y'=\a_ny^n+\cdots+\a_3y^3+\a_2y^2+\a_1y$, when $\gamma'/\gamma=\a_1$ for some non zero $\gamma\in k$ and  $x'\neq \gamma \a_2$ and $x'\neq \gamma^{2}\a_3$ for any $x\in k$. The proof that there are no elements in $E-\ok$ satisfying the differential equation follows immediately from proposition \ref{algsolutions} once we observe that $(y/\gamma)'=$ $\gamma^{n-1}\a_n(y/\gamma)^n+\cdots+\gamma^2\a_3(y/\gamma)^3+\gamma\a_2(y/\gamma)^2$. Let $\C(x)$ be the field of rational functions in one variable over the field of complex numbers with the usual derivation $'=d/dx$. It is evident that $Y=-x$ is a liouvillian solution of
the  differential equation \begin{equation}\label{abeleqn}Y'= \frac{1}{x^3} Y^3+\frac{1}{x^2} Y^2+\frac{1}{x} Y.\end{equation} 
Since $x'/x=1/x$, we see that $y$ is a solution of equation \ref{abeleqn} if and only if $y/x$ is a solution of $Y'=(1/x)(Y^3+Y^2)$. Now it follows from proposition \ref{algsolutions} that equation \ref{abeleqn} has no solution in $E-\overline{\C(x)}$ for any  liouvillian extension  $E$ of  $\overline{\C(x)}$ with $C_E=C_{\overline{\C(x)}}$.


\bibliographystyle{elsarticle-num}
\bibliography{<your-bib-database>}

\end{document}